\documentstyle[11pt,twoside]{article}
\include{graphicx}
\title{On an extension of Watson's lemma due to Ursell}
\author{\sc R. B.\ Paris \\
{\em Division of Computing and Mathematics}, \\
{\em University of Abertay Dundee, Dundee DD1 1HG, UK}
}
\begin{document}
\def\f#1#2{\mbox{${\textstyle \frac{#1}{#2}}$}}
\def\dfrac#1#2{\displaystyle{\frac{#1}{#2}}}
\def\boldal{\mbox{\boldmath $\alpha$}}
{\newcommand{\Sgoth}{S\;\!\!\!\!\!/}
\newcommand{\bee}{\begin{equation}}
\newcommand{\ee}{\end{equation}}
\newcommand{\lam}{\lambda}
\newcommand{\ka}{\kappa}
\newcommand{\al}{\alpha}
\newcommand{\fr}{\frac{1}{2}}
\newcommand{\fs}{\f{1}{2}}
\newcommand{\g}{\Gamma}
\newcommand{\br}{\biggr}
\newcommand{\bl}{\biggl}
\newcommand{\ra}{\rightarrow}
\newcommand{\mbint}{\frac{1}{2\pi i}\int_{c-\infty i}^{c+\infty i}}
\newcommand{\mbcint}{\frac{1}{2\pi i}\int_C}
\newcommand{\mboint}{\frac{1}{2\pi i}\int_{-\infty i}^{\infty i}}
\newcommand{\gtwid}{\raisebox{-.8ex}{\mbox{$\stackrel{\textstyle >}{\sim}$}}}
\newcommand{\ltwid}{\raisebox{-.8ex}{\mbox{$\stackrel{\textstyle <}{\sim}$}}}
\renewcommand{\topfraction}{0.9}
\renewcommand{\bottomfraction}{0.9}
\renewcommand{\textfraction}{0.05}
\newcommand{\mcol}{\multicolumn}
\date{}
\maketitle
\pagestyle{myheadings}
\markboth{\hfill \sc R. B.\ Paris  \hfill}
{\hfill \sc  Watson's lemma\hfill}
\begin{abstract}
In 1991, Ursell gave a strong form of Watson's lemma for the Laplace integral
\[\int_0^\infty e^{-xt}f(t)\,dt\qquad (x\rightarrow+\infty) \]
in which the amplitude function $f(t)$ is regular at the origin and possesses a Maclaurin expansion valid in $0\leq t\leq R$.
He showed that if the asymptotic series for the integral as $x\ra+\infty$ is truncated after $rx$ terms, where $0<r<R$, then the resulting remainder is exponentially small of order $O(e^{-rx})$. In this note we extend this result to include situations when $f(t)$ has a branch point at $t=0$ and when $x$ is a complex variable satisfying $|\arg\,x|<\pi/2$.
\vspace{0.4cm}

\vspace{0.3cm}

\end{abstract}

\vspace{0.3cm}

\noindent $\,$\hrulefill $\,$

\vspace{0.2cm}

\begin{center}
{\bf 1. \  Introduction}
\end{center}
\setcounter{section}{1}
\setcounter{equation}{0}
\renewcommand{\theequation}{\arabic{section}.\arabic{equation}}
Suppose that the function $F(z)$ possesses the asymptotic expansion
\bee\label{e11}
F(z)\sim \sum_{k=0}^\infty c_kz^{-k}
\ee
as $z\ra\infty$ in a sector {\bf S} of the complex $z$-plane. If we denote by $R_n(z)$ the error in approximating $F(z)$ by the sum of the first $n$ terms, then we have exactly
\[F(z)=\sum_{k=0}^{n-1}c_kz^{-k}+R_n(z) \qquad (z\in{\bf S}).\]
By Poincar\'e's definition of an asymptotic series, the error (or remainder) term satisfies the bound $|R_n(z)|\leq A_n|z|^{-n}$, where $A_n$ is an assignable constant for $z$ in the sector {\bf S}. The integer $n$ is arbitrary, but fixed: truncation of the series (\ref{e11}) at the $n$th term consequently yields an error that decays algebraically as $z\ra\infty$ in {\bf S}. When $|z|$ is large, the successive terms in a typical asymptotic expansion (\ref{e11}) initially start to decrease in absolute value, reach a minimum and thereafter increase without bound given the divergent character of the full expansion. If, however, the series is truncated just before this minimum modulus term is reached, then this process is called {\it optimal truncation} and the finite series that results is the {\it optimally truncated expansion}. As we shall see in a specific case below, the resulting rate of decay of $R_n(z)$ is then greatly enhanced.

To illustrate we consider the exponential integral
\[{\cal E}(x):=xe^x\,E_1(x)=\int_0^\infty \frac{e^{-xt}}{1+t}\,dt,\]
where $x>0$, which provides one of the simplest examples of an asymptotic expansion complete with error bound.
Substitution of the identity
\[
\frac{1}{1+t}=\sum_{k=0}^{n-1}(-t)^k+\frac{(-t)^n}{1+t} \qquad (t\neq 1),
\]
followed by termwise integration, yields
\bee\label{e12}
{\cal E}(x)=\sum_{k=0}^{n-1} (-)^k\frac{k!}{x^k}+R_n(x),\qquad R_n(x)=(-)^n\int_0^\infty\frac{t^ne^{-xt}}{1+t}\,dt.
\ee
It is readily seen that 
\[|R_n(x)|<\int_0^\infty  t^ne^{-xt}dt=\frac{n!}{x^n}~,\]
so that this expansion enjoys the property that the remainder is bounded in magnitude by the first neglected term and has the same sign. 

For a given value of $x$, the smallest term in absolute value of the series in (\ref{e12})
occurs when $k=\lfloor x\rfloor$ (except when $x$ is an integer, in which case there are two equally small terms corresponding to $k=x-1$ and $k=x$).
If we denote the optimal truncation index by $N$, then we have
\[{\cal E}(x)=\sum_{k=0}^{N-1} (-)^k\frac{k!}{x^k}+R_N(x).\]
It is clear that $R_N(x)$ is a discontinuous function of $x$, since $N$ changes each time $x$ passes through an integer value. Use of Stirling's formula $N!\sim (2\pi)^{1/2} e^{-N} N^{N+1/2}$ to approximate $N!$ for large $N\simeq x$ then produces the estimate
\[|R_N(x)|<\frac{N!}{x^N}\simeq (2\pi)^\fr\,\frac{e^{-N} N^{N+\fr}}{x^N}\simeq (2\pi x)^\fr e^{-x}\]
as $x\ra+\infty$.
This shows that at optimal truncation the remainder for ${\cal E}(x)$ is of order $x^{1/2}e^{-x}$ as $x\ra+\infty$ and consequently that evaluation of ${\cal E}(x)$ by this scheme will result in an error that is {\it exponentially small\/}. This level of asymptotic approximation has been given the neologism {\it superasymptotics} in \cite{B}; sometimes it is referred to as {\it exponential improvement\/}.

As a rule, it is found that similar exponential improvement at optimal truncation can be achieved in asymptotic expansions of other functions, most notably the confluent hypergeometric functions
which include many well-known special functions, such as Bessel and Airy functions. In these cases the remainder can be written explicitly in the form of an integral and its value estimated to establish its exponentially small nature. The idea of terminating an asymptotic expansion near its smallest term and estimating the remainder is well known in the asymptotic theory of converging factors; see, for example, \cite[\S 7.34]{W} and \cite[Ch. 14]{O}.

There are, however, few general results of this character. In 1958, Jeffreys \cite{HJ} considered
the Laplace transform of $t^\mu f(t)$ (with $\mu>-1,\ f(0)\neq 0$)
which has the expansion for $x\ra+\infty$
\bee\label{e13}
\int_0^\infty e^{-xt} t^\mu f(t)\,dt=\sum_{r=0}^{n-1}u_r+R_n(x),\qquad u_r=\frac{(r+\mu)!}{r! x^{r+\mu+1}}\,f^{(r)}(0),
\ee
where $R_n(x)$ is the remainder after $n$ terms.
If the singularity of $f(t)$ closest to the origin is situated at $Re^{i\alpha}$ ($\alpha\neq 0$), then the least term in modulus in the above series occurs when $n\sim Rx$. When $n$ is so chosen, he showed that the optimal remainder is given by
\[R_n\simeq \frac{u_n}{1-e^{-i\alpha}} \qquad (x\ra+\infty).\]
In the most common case we have $\alpha=\pi$ and the above result shows that the remainder is then approximately half the least term in the expansion.

A significant advance  was made by Ursell \cite{U} who established the following result. 
\newtheorem{theorem}{Theorem}
\begin{theorem}$\!\!\!.$\ \ Let $f(t)$ be analytic in the disc $|t|<R$ with the Maclaurin expansion
\[f(t)=\sum_{n=0}^\infty c_nt^n\qquad(0\leq t\leq R).\]
Let $r$ be a fixed positive quantity such that $0<r<R$ and suppose that $|f(t)|<Ke^{\beta t}$ when $r\leq t<\infty$, where $K$ and $\beta$ are positive constants. Then
\bee\label{e14}
I(x)=\int_0^\infty e^{-xt}f(t)\,dt=\sum_{n=0}^{n_*}\frac{c_n n!}{x^{n+1}}+R_{n_*}(x),
\ee
where, if\/ $n_*$ is chosen to equal $\lfloor rx\rfloor$, the remainder $R_{n_*}(x)=O(e^{-rx})$ as $x\ra+\infty$.
\end{theorem}
Ursell's proof relied on simple bounds for the incomplete gamma functions and was particularly elegant and straightforward; see also \cite[p.~76]{P} for a detailed account.
We observe that the series appearing in (\ref{e13}) and (\ref{e14}), when the upper limit is extended to infinity, is the standard asymptotic expansion of the integrals on the left-hand sides as $x\ra+\infty$ obtained by application of Watson's lemma; see, for example, \cite[p.~44]{DLMF}. Ursell's result shows that if we take $1+\lfloor rx\rfloor$ terms in the asymptotic series, then the resulting approximation to $I(x)$ is exponentially accurate, with the error being $O(e^{-rx})$ as $x\ra+\infty$. 

Our aim in this note is to extend the result in Theorem 1 to situations where the function $f(t)$
possesses a branch point at the origin and to include sectors of the complex $x$-plane. We  carry this out using a modification of the procedure adopted by Ursell.
\vspace{0.6cm}

\begin{center}
{\bf 2. \ The modified expansion}
\end{center}
\setcounter{section}{2}
\setcounter{equation}{0}
\renewcommand{\theequation}{\arabic{section}.\arabic{equation}}
We consider the function $I(z)$ defined by the Laplace integral
\bee\label{e21}
I(z)=\int_0^\infty e^{-zt} f(t)\,dt,
\ee
where $z=xe^{i\theta}$ is a large complex variable with $x>0$ and phase $\theta=\arg\,z$ satisfying $|\theta|\leq\fs\pi-\delta$, $\delta>0$. The amplitude function $f(t)$ satisfies the following conditions:
%\begin{description}
(i) $f(t)$ is holomorphic in the sector {\bf S}: $-\alpha_1+\delta\leq\arg\,t\leq\alpha_2-\delta$;
(ii) $f(t)$ possesses the absolutely convergent expansion
$$f(t)=\displaystyle{\sum_{n=0}^\infty c_n t^{\frac{n+\beta}{\mu}-1}\quad (|t|<R)},$$
where $\mu>0$ and $\Re (\beta)>0$;
and (iii) $|f(t)|\leq Ae^{\sigma |t|} \ \ \ (t\in{\bf S})$,
where $A$ and $\sigma$ are positive quantities.
%\end{description}
Then we have
\begin{theorem}$\!\!\!.$\ \ Let $f(t)$ satisfy conditions (i) -- (iii) above where we suppose that $\alpha_1,\ \alpha_2\geq\fs\pi$. Let $r$ be a fixed positive number such that $r<R$ and $n_*=\lfloor \mu r|z|+\mu-\Re (\beta)\rfloor$. Then
\bee\label{e20}
I(z)=\sum_{n\leq n_*} \frac{c_n}{z^{(n+\beta)/\mu}}\g\bl(\frac{n+\beta}{\mu}\br)+R_{n_*}(z),
\ee
where the remainder $R_{n_*}(z)$ satisfies the bound
\[R_{n_*}(z)=O(e^{-r|z|})\]
as $|z|\ra\infty$ in the sector $|\arg\,z|\leq\fs\pi-\delta$, $\delta>0$.
\end{theorem}

\noindent {\it Proof}.\ \ We rotate the integration path in (\ref{e21}) by $-\theta$ to find
\[I(z)=\int_0^{\infty e^{-i\theta}} e^{-zt}f(t)\,dt=e^{-i\theta}\int_0^\infty e^{-x\tau} f(\tau e^{-i\theta})\,d\tau,\]
where the contribution from the arc at infinity vanishes by assumption (iii) when $|z|>\sigma \mbox{cosec}\, \delta$ \cite[p.~107]{O}. Then, for $0<r<R$,
\begin{eqnarray}
I(z)&=&e^{-i\theta}\sum_{n=0}^\infty c_n\int_0^r e^{-x\tau}(\tau e^{-i\theta})^{(n+\beta)/\mu-1}d\tau+J\nonumber\\
&=&\sum_{n=0}^\infty \frac{c_n}{z^{(n+\beta)/\mu}}\,\gamma\bl(\frac{n+\beta}{\mu}, rx\br)+J\label{e22}
\end{eqnarray}
where $\gamma(a,x)$ is the lower incomplete gamma function and
\[J=e^{-i\theta}\int_r^\infty e^{-x\tau}f(\tau e^{-i\theta})\,d\tau.\]
The sum in (\ref{e22}) converges since 
\[\int_0^r e^{-x\tau} \tau^{\frac{n+\Re (\beta)}{\mu}-1}d\tau<\int_0^r\tau^{\frac{n+\Re (\beta)}{\mu}-1}d\tau=\frac{\mu r^{(n+\Re (\beta))/\mu}}{n+\Re (\beta)}.\]
It may be remarked that (\ref{e22}) represents the first stage of a Hadamard expansion process as discussed in detail in the monograph \cite[Ch. 2]{P}.

Let us now split the series in (\ref{e22}) at $n=n_*$, where $n_*$ is to be specified. Then upon use of the relation
\[\gamma(a,x)+\g(a,x)=\g(a)\]
between the upper and lower incomplete gamma functions, we have
\bee\label{e23}
I(z)=\sum_{n\leq n_*} \frac{c_n}{z^{(n+\beta)/\mu}}\bl\{\g\bl(\frac{n+\beta}{\mu}\br)-\g\bl(\frac{n+\beta}{\mu}, rx\br)\br\}+\sum_{n>n_*}\frac{c_n}{z^{(n+\beta)/\mu}}\gamma\bl(\frac{n+\beta}{\mu}, rx\br)+J.
\ee
If we now define the remainder $R_{n_*}(z)$ by
\[I(z)=\sum_{n\leq n_*} \frac{c_n}{z^{(n+\beta)/\mu}}\g\bl(\frac{n+\beta}{\mu}\br)+R_{n_*}(z),\]
it then follows from (\ref{e23}) that
\[|R_{n_*}(z)|\leq \sum_{n\leq n_*} \frac{|c_n|}{|z^{(n+\beta)/\mu}|}\,\bl|\g\bl(\frac{n+\beta}{\mu}, rx\br)\br|
+\sum_{n>n_*} \frac{|c_n|}{|z^{(n+\beta)/\mu}|}\,\bl|\gamma\bl(\frac{n+\beta}{\mu}, rx\br)\br|+|J|.\]

We now choose
\bee\label{e24}
n_*=\lfloor \mu rx+\mu-\Re (\beta)\rfloor,
\ee
where we recall that $\mu>0$ and $\Re (\beta)>0$ by hypothesis,
and make use of the bounds given in (\ref{a1}) and (\ref{a2})
\begin{eqnarray*}
\bl|\g\bl(\frac{n+\beta}{\mu},rx\br)\br|&\leq & 2e^{-rx} (rx)^{(n+\Re (\beta))/\mu}\qquad (0\leq n\leq n_*)\\
\bl|\gamma\bl(\frac{n+\beta}{\mu},rx\br)\br|&\leq & e^{-rx} (rx)^{(n+\Re (\beta))/\mu}\qquad \ \ (n\geq n_*).
\end{eqnarray*}
Thus we find
\begin{eqnarray}
|R_{n_*}(z)|&\leq& 2\lambda e^{-rx} \sum_{n\leq n_*} |c_n|\, r^{(n+\Re (\beta))/\mu}+\lambda e^{-rx} \sum_{n>n_*}|c_n|\, r^{(n+\Re (\beta))/\mu}+|J|\nonumber\\
&\leq & 2\lambda e^{-rx} \sum_{n=0}^\infty |c_n|\, r^{(n+\Re (\beta))/\mu}+|J|,\label{e25}
\end{eqnarray}
where $\lambda=\exp [\theta \Im (\beta)/\mu]$. By condition (ii) it follows that 
$\sum |c_n| r^{(n+\Re (\beta))/\mu}$ is convergent for $r<R$.

Let $a>\sigma$ be a value of $\Re (z)$ for which $\int_0^\infty e^{-zt} f(t)\,dt$ converges. Then when $x>a$, we find by partial integration \cite[p.~72]{O}
\[J=\int_r^\infty e^{-x\tau} f(\tau e^{-i\theta})\,d\tau=\int_r^\infty e^{-(x-a)\tau}\,e^{-a\tau} f(\tau
 e^{-i\theta})\,d\tau\]
\[=(x-a)\int_r^\infty e^{-(x-a)\tau} {\cal F}(\tau;\theta)\,d\tau,\]
 where
\[{\cal F}(\tau;\theta)=\int_r^\tau e^{-a\tau} f(\tau e^{-i\theta})\,d\tau.\]
If $L(\theta)$ denotes the supremum of $|{\cal F}(\tau;\theta)|$ along the ray $[r,\infty)$ situated in {\bf S} then\footnote{A bound on $|J|$ can also be obtained from condition (iii); see \cite{U}, \cite[p.~76]{P}.}
\bee\label{e26}
|J|\leq (x-a) L(\theta) \int_r^\infty e^{-(x-a)\tau} d\tau=L(\theta)\,e^{-(x-a)r}.
\ee
By combining (\ref{e25}) and (\ref{e26}), we see that with the above choice of $n_*$ the remainder satisfies
$|R_{n_*}(z)|=O(e^{-rx})$, thereby establishing the theorem. \hfill $\Box$

\vspace{0.2cm}

We now suppose that a singularity of $f(t)$ lies on the circle of convergence $|t|=R$ situated at $t_0=Re^{-i\psi}$, where $0<\psi<\fs\pi$. An analogous treatment applies when there is a singularity at $Re^{i\psi}$. During the process of path rotation, it is now possible to cross over the singularity,
thereby receiving either a pole contribution or a branch point contribution depending on the nature of the singularity at $t_0$.
\begin{theorem}$\!\!\!.$\ \ In the case when there is a singularity of $f(t)$ on the circle of convergence $|t|=R$ at the point $t_0=Re^{-i\psi}$, where $0<\psi<\fs\pi$, then the remainder $R_{n_*}(z)$ satisfies the bound
\bee\label{e27}
R_{n_*}(z)=O(e^{-r|z|})+\Upsilon(\theta) \,O(e^{-|z|R \cos (\theta-\psi)}),
\ee
where $n_*$ is defined in (\ref{e24}) when $r<R$ and
\[\Upsilon(\theta)=\left\{\begin{array}{lrl} 0 & -\fs\pi+\delta\leq&\!\!\!\theta\leq\psi\\
\\
[-0.3cm]

1 & \psi+\delta\leq&\!\!\!\theta\leq\fs\pi-\delta.\end{array}\right.\]
An analogous result holds if the singularity is situated in the first quadrant at $t_0=Re^{i\psi}$.
\end{theorem}

\noindent {\it Proof}.\ \ The treatment when $\theta<\psi$ and when $\theta>\psi$ follows that described in the proof of Theorem 2, where in the latter case there is the addition of the contribution from the singularity at $t_0$. In the case of a simple pole this yields a contribution of $O(\exp\,(-zt_0))=O(\exp\,(-|z|R \cos (\theta-\psi)))$. In the case of a branch point singularity, the additional contribution is given by the loop integral
\[e^{-zt_0} \int_\infty^{(0-)} e^{-zu} f(t_0+u)\,du,\]
which is again of $O(\exp\,(-zt_0))$.

When $\theta=\psi$, we can take the integration path for $J$ to be the path commencing at $re^{-i\psi}$ and passing to infinity parallel to the real $t$-axis, thereby avoiding the singularity at $t_0$. Then
\[J=\int_{re^{-i\psi}}^{re^{-i\psi}+\infty}e^{-zt} f(t)\,dt=e^{-rx} \int_0^\infty e^{-zu} f(u+re^{-i\psi})\,du,\]
where on the path $u\in [0, \infty)$ we have by condition (iii) $|f(u+re^{-i\psi})|\leq Ae^{\sigma |u+re^{-i\psi}|}\leq Ae^{\sigma (u+r)}$. Hence
\[|J|\leq Ae^{-rx}\int_0^\infty e^{-\Re (z) u} e^{\sigma (u+r)}du\leq \frac{Ae^{-r(x-\sigma)}}{\Re (z)-\sigma}<A' \frac{e^{-rx}}{x \cos \psi},\]
where $A'$ is a constant. By the argument leading to (\ref{e25}), the remainder satisfies $|R_{n_*}(z)|=O(e^{-rx})$
when $\theta=\psi$ and $n_*$ is chosen according to (\ref{e24}), thereby establishing the theorem. \hfill $\Box$

\vspace{0.2cm}

\noindent {\bf Remark 1}.\ \  If there is a second singularity on $|t|=R$ with $|\arg\,t|<\fs\pi$ then the modification to Theorem 3 is straightforward. Similarly, if there is an additional singularity of $f(t)$  beyond the circle of convergence $|t|=R$ situated at $t_0=R' e^{-i\psi}$, say, where $R'>R$ and $0<\psi<\fs\pi$, then the above argument is easily modified to yield the remainder
\bee\label{e28}
|R_{n_*}(z)|=O(e^{-r|z|})+\Upsilon(\theta) \,O(e^{-|z|R' \cos (\theta-\psi)}).
\ee

\vspace{0.2cm}

\noindent {\bf Remark 2}\,.\ \  The second order estimate in (\ref{e27}) will only be significant (when $\theta>\psi$) if 
\bee\label{e29}
\theta-\psi>\arccos (r/R). 
\ee
A similar remark applies to (\ref{e28}) when $R$ is replaced by $R'$.

\vspace{0.6cm}

\begin{center}
{\bf 3. \ Numerical verification}
\end{center}
\setcounter{section}{3}
\setcounter{equation}{0}
\renewcommand{\theequation}{\arabic{section}.\arabic{equation}}
We present some examples which illustrate the remainder estimates in Theorems 2 and 3 when the underlying asymptotic series is truncated after $n_*$ terms, where $n_*$ is defined in (\ref{e24}).
\bigskip

\noindent Example 1.\ \ We consider the confluent hypergeometric function
\[U(a,a-b+1,z)=\frac{1}{\g(a)}\int_0^\infty e^{-zt} t^{a-1} (1+t)^{-b}dt\qquad (\Re (a)>0)\]
for which 
\[f(t)=t^{a-1}(1+t)^{-b}=\sum_{n=0}^\infty \frac{(-)^n (b)_n}{n!}\,t^{n+a-1} \quad (|t|<1),\]
where $(a)_n=\g(a+n)/\g(a)$ is Pochhammer's symbol.
The function $f(t)$ has a singularity at $t=-1$, so that $R=1$, and is associated with the parameters $\mu=1$, $\beta=a$. Then the remainder is defined by 
\[R_{n_*}(z)=U(a,a-b+1,z)- z^{-a}\sum_{n\leq n_*} \frac{(-)^n (a)_n (b)_n}{n!\,z^{n}},\]
where the truncation index $n_*=r|z|$ with $r<1$.
According to Theorem 2, the remainder $R_{n^*}(z)=O(e^{-r|z|})$ as $|z|\ra\infty$ in $|\arg\,z|\leq\fs\pi-\delta$.

It is known \cite[\S 13.7(iii)]{DLMF} that when $n_*=|z|$ the remainder for the $U$ confluent hypergeometric function satisfies $R_{n_*}(z)=O(e^{-|z|})$ uniformly for $|\arg\,z|\leq\pi$.
\bigskip

\noindent Example 2.\ \ The function
\[I(z)=\int_0^\infty e^{-zt}\,\frac{dt}{\sqrt{1+t^2}}=\fs\pi \,{\bf K}_0(z),\]
where ${\bf K}_0(z)$ is the Struve function defined in \cite[Eq.~(11.2.5)]{DLMF}. The function
\[f(t)=(1+t^2)^{-\fr}=\sum_{n=0}^\infty \frac{(-)^n (\fs)_n}{n!}\,t^{2n} \qquad (|t|<1)\]
has singularities at $t=\pm i$ and is associated with the parameters $\mu=\beta=\fs$. Then the remainder is defined by
\[R_{n_*}(z)=I(z)-\frac{1}{2}\sum_{n\leq n_*} \frac{(-)^n (\fs)_n (\fs)_n}{(\fs z)^{2n+1}},\]
where the truncation index $n_*=r|z|$ with $r<1$. From Theorem 2, $R_{n^*}(z)=O(e^{-r|z|})$ as $|z|\ra\infty$ in $|\arg\,z|\leq\fs\pi-\delta$. 

In Fig.~1 we show the variation of $|R_{n_*}(z)|$ in Examples 1 and 2 as a function of $\theta$ in the range $0\leq\theta<\fs\pi$
for different values of $x$. It can be seen that each curve lies close to the value of the order estimate $e^{-r|z|}$ and that (on the scale of the figure) there is little variation in the value of $|R_{n_*}(z)|$ over
the range $0\leq\theta<\fs\pi$. Conjugate values are obtained for $-\fs\pi<\theta\leq 0$.
\begin{figure}[t]
	\begin{center}	{\tiny($a$)}\includegraphics[width=0.4\textwidth]{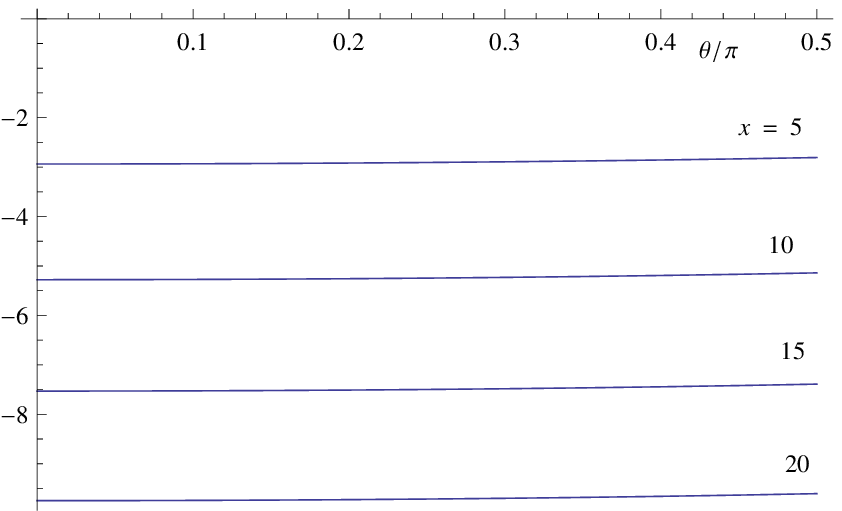}\qquad
	{\tiny($b$)}\includegraphics[width=0.4\textwidth]{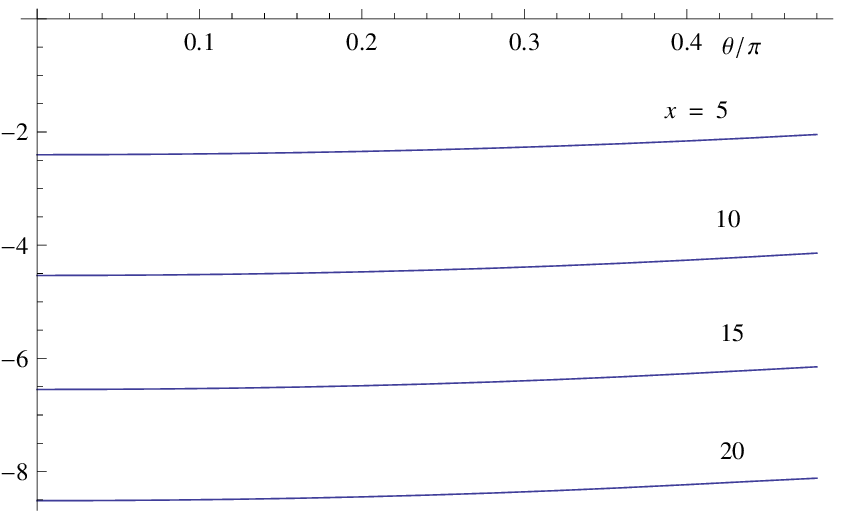}\\
\caption{\small{Plots of $|R_{n_*}(z)|$ against $\theta/\pi$ (on a $\log_{10}$ scale) for $x=5, 10, 15, 20$
and $r=0.8$ for (a) the confluent hypergeometric function $U(a, a-b+1,z)$ with $a=\fs$, $b=\frac{3}{4}$ and (b)
the Struve function in Example 2.}}
	\end{center}
\end{figure}
\bigskip

\noindent Example 3.\ \ The function
\[I(z)=\int_0^\infty \frac{e^{-zt}dt}{e^{-i\psi}-t}dt\qquad (0<\psi<\fs\pi)\]
is associated with the amplitude function
\[f(t)=e^{i\psi} (1-e^{i\psi}t)^{-1}=\sum_{n=0}^\infty e^{i\psi (n+1)} t^n\qquad (|t|<1),\]
so that $\mu=\beta=1$. The function $f(t)$ has a simple pole at $t_0=e^{-i\psi}$, which produces  the residue contribution $2\pi i \exp\,(-ze^{-i\psi})$. Then, with the remainder defined by
\[R_{n_*}(z)=I(z)-\sum_{n\leq n_*} n! (ze^{-i\psi})^{-n-1},\]
we have from Theorem 3 that as $|z|\ra\infty$  
\bee\label{e31}
R_{n_*}(z)=O(e^{-rx})+\Upsilon(\theta)\,O(e^{-|z| \cos (\theta-\psi)}).
\ee

\bigskip

\noindent Example 4.\ \ Our final example is the function
\[I(z)=\int_0^\infty \frac{e^{-zt}dt}{\sqrt{1-te^{i\psi}}}\qquad (0<\psi<\fs\pi).\]
Here we have
\[f(t)=(1-te^{i\psi})^{-\fr}=\sum_{n=0}^\infty \frac{(\fs)_n}{n!}\, (te^{i\psi})^n\qquad (|t|<1),\]
so that $\mu=\beta=1$. The amplitude function $f(t)$ possesses a square-root branch point at $t_0=e^{-i\psi}$, which produces the contribution
\[\int_\infty^{(t_0-)}\frac{e^{-zt}dt}{\sqrt{1-te^{i\psi}}}=e^{-zt_0}\int_{-\infty}^{(0+)} \frac{e^{zu}du}{(ue^{i\psi})^\fr}=2ie^{-zt_0} \sqrt{\frac{\pi}{ze^{i\psi}}}\]
during the path rotation argument. With the remainder defined by
\[R_{n_*}(z)=I(z)-\sum_{n\leq n_*} (\fs)_n e^{in\psi} z^{-n-1},\]
we have the order estimate\footnote{We omit the $|z|^{-1/2}$ term in the second exponential.} given in (\ref{e31}).

In Fig.~2(a) we give an example of the variation of $\log_{10}(e^{r|z|}\,|R_{n_*}(z)|)$ for the integral in Example 3 as a function of $\theta$ when $x=20$ and $\psi=0.4\pi$. For this value of $\psi$, the condition (\ref{e29}) is not satisfied and so the behaviour of $|R_{n_*}(z)|$ is controlled by the $O(e^{-r|z|})$ term. When $\psi=0.1\pi$, however, the condition (\ref{e29}) yields $\theta-\psi>\arccos\,0.8 \simeq 0.205\pi$. Thus, for $\theta\, \gtwid\,0.305\pi$
the behaviour of $|R_{n_*}(z)|$ is controlled by the pole contribution, as can be seen from the figure where the value changes by about five orders of magnitude. The dashed curve shows the variation of $e^{|z| \cos (\theta-\psi)}\,|R_{n_*}(z)|$: it is seen that this curve approaches the value $\log_{10}2\pi\doteq 0.79818$ 
confirming the dominance of the second term in (\ref{e31}) for $\theta\,\gtwid\,0.3\pi$.
Fig.~2(b) shows the variation of the remainder in Example 4 for the same values of $r$, $x$ and $\psi$ as in Example 3.
Similar remarks apply in this case.

\begin{figure}[t]
	\begin{center}	{\tiny($a$)}\,\,\includegraphics[width=0.4\textwidth]{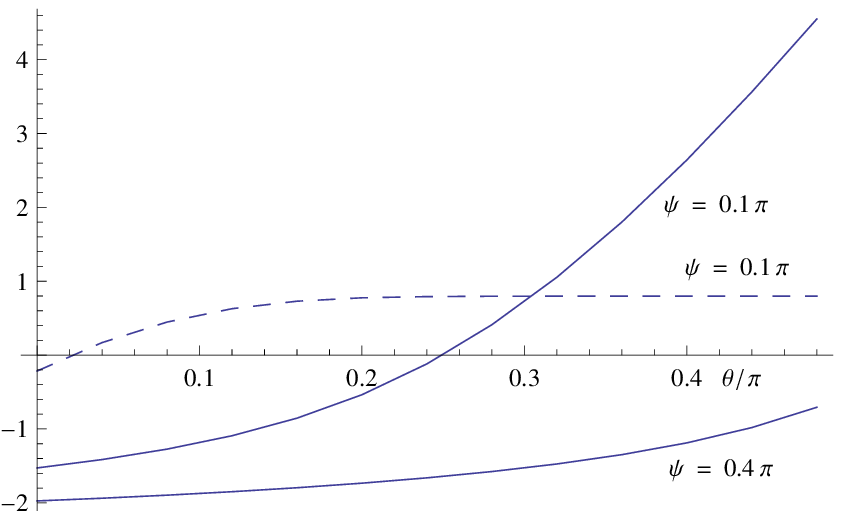}\qquad
	{\tiny($b$)}\includegraphics[width=0.4\textwidth]{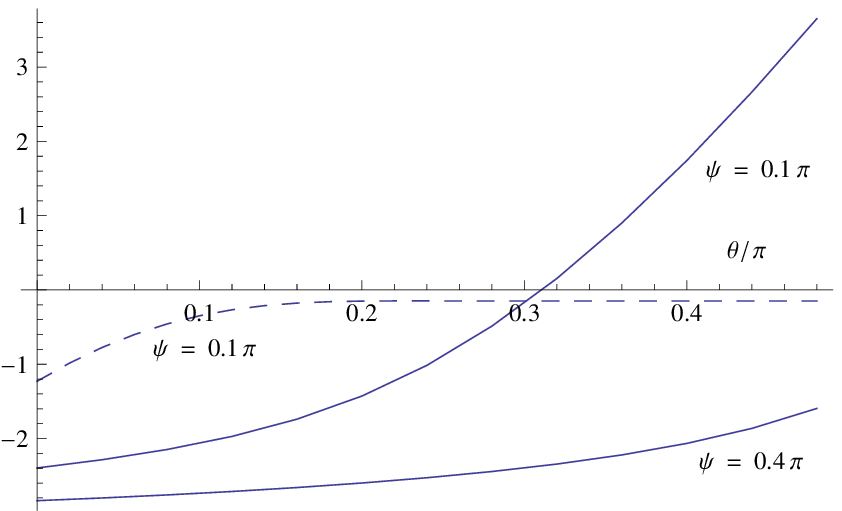}\\
\caption{\small{Plots of $\log_{10}(e^{r |z|}|R_{n_*}(z)|)$ against $\theta/\pi$ (solid curves) for $r=0.8$
and two values of $\psi$ when $x=20$: (a) for Example 3 and (b) for Example 4. The dashed curves show the variation of $\log_{10}(e^{|z|\cos (\theta-\psi)} \,|R_{n_*}(z)|)$ as a function of $\theta/\pi$ when $\psi=0.1\pi$.}}
	\end{center}
\end{figure}
\vspace{0.6cm}

\begin{center}
{\bf Appendix: \ Derivation of bounds on the incomplete gamma functions}
\end{center}
\setcounter{section}{1}
\setcounter{equation}{0}
\renewcommand{\theequation}{\Alph{section}.\arabic{equation}}
Let $\omega=a+ib$, where $a$ and $b$ are real, and suppose that $\chi>0$. Then, following Ursell \cite{U}, we have
\begin{eqnarray*}
|\g(\omega+1,\chi)|&=&\bl|\int_\chi^\infty e^{-t}t^\omega dt\br|=e^{-\chi} \chi^{a+1} \bl|\int_0^\infty e^{-u\chi} (1+u)^\omega du\bl|\\
&=& e^{-\chi} \chi^{a+1} \int_0^\infty e^{-u\chi}(1+u)^a du.
\end{eqnarray*}
When $-1\leq a\leq 0$, we have
\[|\g(\omega+1,\chi)|\leq e^{-\chi} \chi^{a+1}\int_0^\infty e^{-u\chi}du\leq e^{\chi} \chi^{a+1} \quad (\chi\geq 1).\]
When $0\leq a\leq 1$, we have $(1+u)^a\leq 1+u$ for $u\geq 0$ so that
\begin{eqnarray*}
|\g(\omega+1,\chi)|&\leq & e^{-\chi} \chi^{a+1}\int_0^\infty e^{-u\chi}(1+u)\,du\leq e^{\chi} \chi^{a+1}\bl(\frac{1}{\chi}+\frac{1}{\chi^2} \br)\\
&\leq & 2e^{-\chi} \chi^{a+1}\qquad (\chi\geq 1).
\end{eqnarray*}
When $a\geq 1$, we have
\begin{eqnarray*}
|\g(\omega+1,\chi)|&\leq & e^{-\chi} \chi^{a+1}\int_0^\infty \{e^{-u}(1+u)\}^a\,du\\
&\leq& e^{-\chi} \chi^{a+1} \int_0^\infty e^{-u}(1+u) du=2e^{-\chi} \chi^{a+1},
\end{eqnarray*}
the first inequality holding when $a\leq\chi$ and the second inequality when $a\geq 1$, since the integrand $e^{-u}(1+u)\leq 1$ on $[0,\infty)$. The resulting bound therefore holds for $1\leq a\leq\chi$.
Collecting together these results, we therefore obtain the upper bound
\bee\label{a1}
|\g(\omega+1,\chi)|\leq 2e^{-\chi} \chi^{a+1} \qquad (-1\leq a\leq \chi)
\ee
provided $\chi\geq 1$.

Proceeding in a similar manner for the lower incomplete gamma function, we have
\begin{eqnarray*}
|\gamma(\omega+1,\chi)|&=& \bl|\int_0^\chi e^{-t}t^\omega dt\br|=e^\chi \chi^{a+1} \bl|\int_0^1 e^{u\chi} (1-u)^\omega du\br|\\
&\leq& e^{-\chi} \chi^{a+1} \int_0^1 e^{u\chi} (1-u)^adu.
\end{eqnarray*}
It is easily seen that $0\leq e^{u\chi}(1-u)^a\leq 1$ on $[0,1]$ when $a\geq \chi$, and hence we obtain the bound
\bee\label{a2}
|\gamma(\omega+1,\chi)|\leq e^{-\chi} \chi^{a+1} \qquad(a\geq\chi).
\ee

\end{document}